\documentclass[10pt]{book}
\usepackage{array}
\makeatletter
\newif\if@borderstar
\def\bordermatrix{\@ifnextchar*{%
\@borderstartrue\@bordermatrix@i}{\@borderstarfalse\@bordermatrix@i*}%
}
\def\@bordermatrix@i*{\@ifnextchar[{\@bordermatrix@ii}{\@bordermatrix@ii[() ]}}
\def\@bordermatrix@ii[#1]#2{%
\begingroup
\m@th\@tempdima8.75\p@\setbox\z@\vbox{%
\def\cr{\crcr\noalign{\kern 2\p@\global\let\cr\endline }}%
\ialign {$##$\hfil\kern 2\p@\kern\@tempdima & \thinspace %
\hfil $##$\hfil && \quad\hfil $##$\hfil\crcr\omit\strut %
\hfil\crcr\noalign{\kern -\baselineskip}#2\crcr\omit %
\strut\cr}}%
\setbox\tw@\vbox{\unvcopy\z@\global\setbox\@ne\lastbox}%
\setbox\tw@\hbox{\unhbox\@ne\unskip\global\setbox\@ne\lastbox}%
\setbox\tw@\hbox{%
$\kern\wd\@ne\kern -\@tempdima\left\@firstoftwo#1%
\if@borderstar\kern2pt\else\kern -\wd\@ne\fi%
\global\setbox\@ne\vbox{\box\@ne\if@borderstar\else\kern 2\p@\fi}%
\vcenter{\if@borderstar\else\kern -\ht\@ne\fi%
\unvbox\z@\kern-\if@borderstar2\fi\baselineskip}%
\if@borderstar\kern-2\@tempdima\kern2\p@\else\,\fi\right\@secondoftwo#1 $%
}\null \;\vbox{\kern\ht\@ne\box\tw@}%
\endgroup}
\makeatother

\usepackage{amsmath,amssymb,latexsym}
\newtheorem{thm}{Theorem}

\newtheorem{lem}[thm]{Lemma}

\newcommand{\be}{\begin{eqnarray*}}
\newcommand{\ee}{\end{eqnarray*}}
\newcommand{\ben}{\begin{eqnarray}}
\newcommand{\een}{\end{eqnarray}}

\newcommand{\Rex}{\mathbb R}

\newcommand\I{\mathcal{I}}
\newcommand\Q{\mathcal{Q}}

\newcommand{\teij}{\theta_{ij}}

\DeclareMathOperator{\Di}{Di}

\DeclareMathOperator{\BF}{BF}

\begin{document}
\renewcommand{\baselinestretch}{1.2}
\markright{
}
\markboth{\hfill{\footnotesize\rm G. Consonni and G. Pistone
}\hfill}
{\hfill {\footnotesize\rm Bayesian analysis with zero-probability cells} \hfill}
\renewcommand{\thefootnote}{}
$\ $\par
\fontsize{10.95}{14pt plus.8pt minus .6pt}\selectfont
\vspace{0.8pc}
\centerline{\large\bf Algebraic Bayesian analysis of contingency tables}
\vspace{2pt}
\centerline{\large\bf with possibly zero-probability cells}
\vspace{.4cm}
\centerline{Guido Consonni and Giovanni Pistone}
\vspace{.4cm}
\centerline{\it University of Pavia, Italy and Politecnico di Torino, Italy}
\vspace{.55cm}
\fontsize{9}{11.5pt plus.8pt minus .6pt}\selectfont

\begin{quotation}
\noindent {\it Abstract:}
  In this paper we consider a Bayesian analysis of contingency tables
  allowing for the possibility that cells may have probability
  zero. In this sense we depart from standard
  log-linear modeling that implicitly assumes a positivity constraint.
  Our approach leads us to consider mixture models for contingency tables,
  where the components of the mixture, which we call model-instances,
  have distinct support.
  We rely on ideas from polynomial algebra in order to identify the
  various model instances. We also provide a method  to
  assign prior probabilities to each instance of the model, as well as
  describing methods for constructing priors on the parameter space of
  each instance.  We illustrate our methodology through
  a $5 \times 2$ table involving two structural zeros, as well as a zero count.
  The results we obtain show that our analysis may lead to conclusions that are substantively different
  from those that would obtain in a standard framework, wherein   the possibility of zero-probability cells
 is not explicitly accounted for.\par

\vspace{9pt}
\par
\end{quotation}\par

\fontsize{10.95}{14pt plus.8pt minus .6pt}\selectfont

\setcounter{chapter}{1} \label{sec:intro}
\setcounter{equation}{0} 
\noindent {\bf 1. Introduction}

The analysis of contingency tables has a well established tradition,
both in the frequentist and Bayesian setting.  A typical framework
for this analysis is represented by the exponential family
representation of the sampling distribution, together with the
log-linear, or more generally log-affine, model for the expected
cell count, see Lauritzen (1996, ch. 4) for a rigorous treatment.
Under multinomial sampling, this approach presupposes implicitly
that  cell-probabilities, equivalently  cell-expected counts, are
strictly positive. On the other hand, this assumption is not
particularly justified from a substantive viewpoint; indeed, as we
shall argue below, it might well hide some interesting aspects of
modeling.

Typically, the positivity constraint is viewed as problematic
when performing Maximum Likelihood Estimation (MLE) in a log-linear
framework if there are some cells  having zero counts, see for instance the
discussion in Christensen (1997, ch. 8).  One usually distinguishes
between \emph{structural} (or ``fixed'') zeros, and \emph{random} (or
``sampling'') zeros. The former arise when the cells are logically
forced to have a zero-count.  Consider for instance a
cross-classification for people where the personal highest educational
attainment (Less than high school, High school, College, Postgraduate)
is recorded at a given time, and five years later. Clearly it is
impossible for someone to have a highest attainment of College on the
first time point, and Less than high school or High school five years
later; in general every cell that corresponds to lower attainment at
the second time period compared to the first time period is a
structural zero. On the other hand, random zeros are
typically thought to occur either because the sample size, or the
corresponding cell probability, or both are ``small'', as
it occurs in sparse contingency tables.

Structural zeros are typically dealt with by removing them altogether
from the analysis.  One way to do this is through regression models on
effect codings, see e.g. Simonoff (2003, sect. 6.4).  Random zeros on
the other hand require special handling. Essentially one should first
identify those cells for which the regular MLE
of the cell-probability does not exist, i.e. is zero (this requires
special care as such cells need not coincide with those having zero
counts), and then remove them from the analysis.  In any case the
computation of the degrees of freedom for model testing must be done
on a case by case basis, and requires some ingenuity.  Another
difficulty  generated by the presence of random zeros  is that asymptotic arguments
may effectively break down because of the small-sample size, although
some computer programs may still provide MLEs when they actually do
not exist. For an informative account of the above problems see
Haberman (1974),
Bishop et al. (1975, sect. 5) and
Christensen (1997, sect. 8.3).  Recently Eriksson et al. (2006)
have provided a polyhedral description of the conditions for the
existence of the MLE for a hierarchical log-linear model together with an algorithm
for determining if the MLE exists.

In this paper we take the view that modeling of contingency tables
should allow explicitly for the possibility of zero-probability cells
not only to deal with structural zeros but also with zero-counts whose
nature is undecided, in the sense that their occurrence may be
consistent with either a zero probability or a positive probability:
we call these cells \emph{possibly zero-probability} cells.

An early paper that takes a similar view is Lauritzen (1975), although the techniques
used there are quite different from the ones that we employ here.

From a modeling perspective, we contend that, for each
given model, the usual exponential-family/log-linear representation of
the sampling distribution is simply one \emph{instance} of such model, while
several other instances are conceptually consistent with the assumed
model, each being essentially  a log-linear model with a restricted support.
The identification of such instances  represent a crucial aspect in the
implementation process, and is typically of high complexity.

In our work we rely on ideas from
polynomial algebra and the related geometric and combinatorial structure, which  have
been recently applied to the analysis of some classes of (finitely) discrete
statistical models. In particular,  Eriksson et al. (2006) deal
with hierarchical log-linear models, while  Geiger et al. (2006)
discuss graphical models.

Our approach falls broadly under the heading of
\emph{Algebraic Statistics}, see Pistone et al. (2001)
for an early general account, as well as the pioneering work of Diaconis and
Sturmfels (1998).
The field is now growing at an impressive speed both
in terms of theoretical contributions and applications, see for
example the recent monograph by Pachter and Sturmfels (2005). Further useful references are
Geiger et al. (2001), who develop the concept of stratified
exponential families,
as well as Garcia, Stillman and
Sturmfels (2005) who carry out the analysis of Bayesian networks
from an algebraic statistical perspective.
Rapallo
(2006) discusses some basic algebraic statistics tools that deal
explicitly with models for contingency tables and represents a simple
and useful introduction to this paper.
Our interest in the use  of algebraic
methodology for statistical purposes was stimulated by the  availability of various symbolic
computational software: here we use CoCoA developed and maintained at
the University of Genova, Italy. An other option could be the softare 4ti2.

A specific feature of this paper is the combination of methods from
algebraic statistics with the Bayesian approach. Specifically, we shall
deal with issues like the assignment of a prior on model space, prior
elicitation on the parameter space under each model, or instances thereof;
together with model choice using the Bayes factor, see
Kass and Raftery (1995) for a review.

The paper is organized as follows: Section 2 contains some basic tools
from algebraic statistics that are used in the paper; in Section 3 such tools are applied to a real data-set; Section 4 is the core of the paper, presenting a Bayesian approach testing  quasi-independence in two-way contingency tables using a mixture of model-instances, thus accounting for the possible presence of zero-probability cells. Finally, Section 5
 summarizes the paper and presents some points for discussion.

\setcounter{chapter}{2} 
\setcounter{equation}{0} 
\noindent {\bf 2. Algebraic statistical models} 

Consider a finite state space $\Q$ and a probability distribution on
$\Q$, which we can write as $\{ p(x), \, x \in \Q  \}$,
with $p(x) \geq 0$  and $\sum_{x \in \Q} p(x) = 1$.
In particular, we shall deal with multi-way contingency tables identified by a
collection of factors $X=\{X_1,\ldots,X_F  \}$. If   $\I_f$ denotes the
set of levels for the factor $X_f$, $f=1,\ldots,F$,
the state space is a product space, i.e  $\Q= \times_{f=1}^F \I_f$.

 A log-linear
model assumes that $p(x)>0$ and that $\log p(x)$ belongs to a
linear subspace $H$ of $L=\Rex^{\Q}$, where $\Rex^{\Q}$ denotes as
usual the vector space of real-valued functions on $\Q$. If $H$ is spanned by
$\{T_1, \ldots, T_s\}$, where the $T_j$'s are integer valued
functions,
we can write the log-linear model as
\ben
\label{formula:log linear}
\log p(x)=\sum_{j=1}^s(\log\zeta_j)T_j(x),
\een
with $\sum_x p(x)=1$.
%
Recall that \eqref{formula:log linear}
assumes strict positivity of $p(x)$.
 However the latter is no longer needed if we rewrite \eqref{formula:log linear} as
%
%
\begin{equation} \label{formula:toric-2}
  q(x) = \zeta_1^{T_1(x)}\cdots \zeta_s^{T_s(x)}, \quad \zeta_j \geq 0, \quad j=0,\ldots,s,
\end{equation}
where $q(x)$ is the un-normalized probability, so that the parameters $\zeta_1,\dots, \zeta_s$ are
only subject to non-negativity constraints.
Notice that \eqref{formula:toric-2} is,
for each $x \in \Q$, a (monic) monomial in the indeterminates  $\zeta_1,\dots, \zeta_s$.
 When $x$ scans
$\Q$, we get a system of binomial equations and so
\eqref{formula:toric-2} could also be called a parametric
\emph{toric} model, borrowing terminology from commutative algebra,
see Sturmfels (1996), as suggested in Pistone et al. (2001).

When the cell probabilities are assumed to be strictly positive, then
the log-linear model \eqref{formula:log linear} and  the toric
model \eqref{formula:toric-2} can be easily shown to be
equivalent. A third expression of the same
model can be derived by elimination of the indeterminates $\zeta_1,\dots,\zeta_s$ in the monomial
parameterization of equation \eqref{formula:toric-2}. In fact, if
$
M =
\begin{bmatrix}
  T_1(x) & \cdots & T_s(x)
\end{bmatrix}_{x \in \Q}
$
is the design matrix of the log-linear model of equation
\eqref{formula:log linear},
 the orthogonal space of its range can be generated by integer valued vectors with zero sum $K =
 \begin{bmatrix}
   k_1 & \cdots & k_r
 \end{bmatrix}
 $, and equation \eqref{formula:toric-2} gives for each $j=1,\dots,r$
\begin{equation} \label{formula:toric-3}
  \prod_{x} q(x)^{k_j(x)} = \prod_{x}
  \left(\zeta_1^{T_1(x)}\cdots \zeta_s^{T_s(x)}\right)^{k_j(x)}=
  \zeta_1^{T_1(x)\cdot k_j(x)}\cdots \zeta_s^{T_s(x) \cdot k_j(x)} = 1,
\end{equation}
where the dot symbol \lq \lq $\cdot$\rq \rq denotes scalar product.

As the sum of the elements of each $k_j$, $j = 1,\dots,r$, is zero,
the sum of the elements of both the positive part $k_j^+$ and the negative part $k_j^-$ are equal,
so that we could write equation \eqref{formula:toric-3} as
\begin{equation} \label{formula:implicit}
  \prod_{x} q(x)^{k_j(x)^+} - \prod_{x} q(x)^{k_j(x)^-} = 0,
  \qquad j=1, \ldots, r.
\end{equation}

It follows that the toric model \eqref{formula:toric-2} implies a set of $r$ binomial and
homogeneous equations in the un-normalized probabilities $q(x)$, $x
\in \mathcal Q$.

If the probabilities are assumed to be strictly positive, then the
three descriptions, i.e. log-linear \eqref{formula:log linear}, toric
\eqref{formula:toric-2} and implicit binomial
\eqref{formula:implicit}, are equivalent. We remark that while \eqref{formula:log linear} and \eqref{formula:toric-2} are parametric models,
the nature of \eqref{formula:implicit} is essentially non-parametric.
When the positivity assumption is relaxed, a non trivial situation
occurs. The basic fact is that different toric parameterizations can
lead to the same implicit binomial, because they are equivalent only
on the strictly positive part of the model. However,   the implicit
binomial equations are satisfied by all  limits of the positive
cases; thus the implicit binomial is the best expression of the so
called extended exponential model, i.e. the exponential model plus
all its limits.

We summarize here  a few basic facts of the theory of toric
statistical models. Given a log-linear model and all its limit points,
a specific set of configurations of zero-probability cells arises. This
set cannot be recovered by setting to zero some parameters
in a generic toric parametric representation, because most
of the equivalent toric representations will not produce all possible
probabilities of the model in Equation \eqref{formula:implicit}. However, there exists a
``maximal'' parametric toric representation, such that all configurations
of zero-probability cells compatible with, i.e. limit of,  the initial model
are obtained by letting some parameters be zero. Such
representation results from the following steps:
\begin{enumerate}
\item All toric models compatible with the implicit binomial model
  \eqref{formula:implicit} are
  characterized by a string of  $T$'s exponents, see \eqref{formula:toric-2}, which is a non-negative
  integer vector orthogonal to the basis $[k_1 \dots k_r]$ of the
  orthogonal space of the initial design matrix $M$.
\item The lattice of non-negative integer vectors
$t \in \mathbb N_+^\Q$
 such that the condition
   $t \cdot k_j = 0$ holds for each  $j=1,\dots,r$,  has a finite number of
  generators that can be computed with symbolic software. Here ``generator'' means that all such vectors are
  component-wise sums of a finite number of generators, possibly
  repeated. The minimal set of generators is called minimal Hilbert basis.

\item If the generators are $S_1, \dots , S_u$, then the ``maximal'' toric model is
  \begin{equation}
    \label{formula:maxtoric}
    q(x) = \zeta_1^{S_1(x)}\cdots \zeta_u^{S_u(x)} \quad x \in \Q.
  \end{equation}
\end{enumerate}
Here ``maximal'' means that \eqref{formula:maxtoric} is a (possibly non-identifiable) parameterization
of the full implicit binomial model, i.e. the extended model. All
members of the implicit model \eqref{formula:implicit} with zero-cell probabilities are obtained by letting some
$\zeta_j$'s  be zero. Assume e.g. we let $\zeta_1 = 0$.
Then the support of the resulting probability will be the set $\Q_1 = \{x \in \Q : S_1(x) =0 \}$.
On such a restricted support, the model will be again toric:
\begin{equation*}
    q(x) = \zeta_2^{S_2(x)}\cdots \zeta_u^{S_u(x)} \quad x \in \Q_1.
\end{equation*}
or exponential if all the other parameters $\zeta_2, \cdots,\zeta_u$
are assumed to be strictly positive. In this sense, we say that each toric model is a union of exponential models with different supports. Each one of these models is called an \emph{instance} of the model.%

Current symbolic software allows to compute, for a given
parametric model, the set of corresponding implicit binomial descriptions.
Moreover, the collection of allowable models obtained by setting some cell
probabilities equal to zero can be identified in terms of the
functions $T_j(x)$, see Geiger et al. (2006) and Rapallo (2006).

\setcounter{chapter}{3} \label{sec:cancer}  
\setcounter{equation}{0} 
\noindent {\bf 3. Example:  new cancer incidence and gender}

We now turn to the discussion of a real example involving both
structural and random zeros. Our analysis aims primarily at illustrating  the main features of our method.

The Division of Cancer Prevention and Control of the National Cancer
Institute in the United States provides (estimates of) counts of new cases of cancer
classified according to various demographic and geographic factors,
see Simonoff (2003, p. 226).  The following table reports data for
different types of cancer separated by gender for Alaska in year
1989.

\begin{center}
\begin{tabular}{l|rr|r}
\emph{Type of cancer}& Female & Male & Total \\ \hline
Lung & 38 & 90 & 128\\
Melanoma & 15 & 15 & 30 \\
Ovarian & 18 & * & 18 \\
Prostate & * & 111 & 111 \\
Stomach & 0 & 5 & 5 \\ \hline
Total & 71 & 221 & 292
\end{tabular}
\end{center}
Clearly cells $(3,2)$ and $(4,1)$ are structural zeros, while we
regard the zero count corresponding to the combination (Stomach,
Female) as a possibly zero-probability cell. A typical assumption that
is of interest in this case is that of \emph{quasi-independence} ($QI$),
corresponding to the standard independence assumption for all cells, excluding
those having a  structural zero.
For this hypothesis, Simonoff (2003, p. 228) finds a $p$-value between 2\% and 3\%,
depending on the  method that is employed. Using a conventional frequentist interpretation,
 the data thus seem to provide  significant evidence against the $QI$-model, although this evidence is not
 very strong.

 Let $I=\{1,2,3,4,5\}$,
$J=\{1,2\}$ denote the set of levels for the rows and columns
respectively, and consider the two-way table with cells in the set $A
= I\times J \setminus \{(3,2),(4,1)\}$, i.e. with cells $(3,2)$ and $(4,1)$
missing.

Under the $QI$-model the
un-normalized cell probabilities $q_{ij}$ are given by
\begin{equation} \label{eq:model1}
  q_{ij} = \rho_i\psi_j, \quad (i,j) \in A.
\end{equation}

If the probabilities are strictly positive, one can take the logarithm of (\ref{eq:model1}), obtaining
\begin{equation*}
\log q_{i,j} = \alpha_i + \beta_j, \quad (i,j) \in A
\end{equation*}
with $\alpha_i = \log \rho_i$, $\beta_j=\log \psi_j$. Accordingly the  design  matrix $M$, together with a suitable choice of an orthogonal matrix $K$, as described in Step 1 of Section 2, are
\begin{equation*}
M =
\bordermatrix[{[]}]{%
& \alpha_1 & \alpha_2 & \alpha_3 & \alpha_4 & \alpha_5 & \beta_1 & \beta_2 \cr
11 & 1 & 0 & 0 & 0 & 0 & 1 & 0 \cr
21 & 0 & 1 & 0 & 0 & 0 & 1 & 0 \cr
31 & 0 & 0 & 1 & 0 & 0 & 1 & 0 \cr
51 & 0 & 0 & 0 & 0 & 1 & 1 & 0 \cr
12 & 1 & 0 & 0 & 0 & 0 & 0 & 1 \cr
22 & 0 & 1 & 0 & 0 & 0 & 0 & 1 \cr
42 & 0 & 0 & 0 & 1 & 0 & 0 & 1 \cr
52 & 0 & 0 & 0 & 0 & 1 & 0 & 1
}%
\qquad
K =
\bordermatrix[{[]}]{%
 & k_1 & k_2 \cr
11 &  1 & 0 \cr
21 & -1 &-1 \cr
31 &  0 & 0 \cr
51 &  0 & 1 \cr
12 & -1 & 0 \cr
22 &  1 & 1 \cr
42 &  0 & 0 \cr
52 &  0 &-1
}%
\end{equation*}

One can check that, under the condition $q_{ij} > 0$, $(i,j) \in A$,
the
model of quasi-independence in \eqref{eq:model1} is equivalent to the
implicit binomial model given by the two constraints
\begin{equation} \label{formula:implicit-binom}
  \left\{ \begin{aligned}q_{11}q_{22} - q_{21}q_{12} &= 0\\
q_{51}q_{22} - q_{21}q_{52} &= 0.
\end{aligned} \right.
\end{equation}
The above equations are the standard conditions for independence in
the two $2 \times 2$ tables with rows $\{1,2\}$, respectively  $\{2,5\}$. This
is equivalent to the independence of the sub-table
$\{1,2,5\}\times\{1,2\}$, since independence for an $R \times
C$-table is equivalent to the  its $2\times2$
minors being zero.

The maximal design matrix $M_{\max}$ and the model in monomial form, see \eqref{formula:maxtoric}, are
\begin{equation} \label{eq:maximal}
M_{\max} =
\bordermatrix[{[]}]{%
& \zeta_1 & \zeta_2 & \zeta_3 & \zeta_4 & \zeta_5 & \zeta_6 & \zeta_7 \cr
11 & 0 & 0 & 0 & 0 & 1 & 0 & 1 \cr
21 & 0 & 0 & 1 & 0 & 0 & 0 & 1 \cr
31 & 1 & 0 & 0 & 0 & 0 & 0 & 0 \cr
51 & 0 & 0 & 0 & 1 & 0 & 0 & 1 \cr
12 & 0 & 0 & 0 & 0 & 1 & 1 & 0 \cr
22 & 0 & 0 & 1 & 0 & 0 & 1 & 0 \cr
42 & 0 & 1 & 0 & 0 & 0 & 0 & 0 \cr
52 & 0 & 0 & 0 & 1 & 0 & 1 & 0
}
\qquad
    \left\{\begin{aligned}
    q_{11} &= \zeta_5 \zeta_7   \\
    q_{21} &= \zeta_3 \zeta_7   \\
    q_{31} &= \zeta_1   \\
    q_{51} &= \zeta_4 \zeta_7 \\
    q_{12} &= \zeta_5 \zeta_6   \\
    q_{22} &= \zeta_3 \zeta_6  \\
    q_{42} &= \zeta_2   \\
    q_{52} &= \zeta_4  \zeta_6
  \end{aligned} \right.
\end{equation}
Notice that the cells associated to a structural zero in the same row
are parameterized independently from the rest of the table. If we take
out these cells, we simply get
 the full independence model on the sub-table with rows
$\{ 1,2,5 \}$.

The instances for the  $QI$-model are computed by considering the
$(2^3-1)(2^2-1)=21$ instances corresponding to independence in the $3\times2$ sub-table,
times the $2^2=4$ instances of the two free cells, plus the $(2^2-1)$
instances where the $3\times2$ sub-table is zero. The total is
87.

\setcounter{chapter}{4} 
\setcounter{equation}{0} 
\noindent {\bf 4. Testing quasi-independence in the new cancer data}

We provide a Bayesian analysis of these data using the methodology developed in the
previous sections. We refer to the model which imposes no restriction on the cell probabilities, save the zero-probability cells $(3,2)$ and $(4,1)$, as the Structural
Zero model and  label it with the symbol $SZ$. Since the table has 10
probability cells, of which 2 are fixed to be zero, the number of
$SZ$-instances is equal to $2^8-1=255$ corresponding to all  possible
combinations of \lq \lq$+$\rq \rq{} and \lq \lq$0$\rq \rq{} in the 8 free cells, excluding the
trivially impossible case of all \lq \lq$0$\rq \rq{}.

Moreover, only two of the above $SZ$-instances
are logically consistent with the observed data: that
giving a positive probability to all  eight free cells; and that
giving zero-probability to cell $(5,1)$ only.  We label these instances
$SZ_{0}$ and $SZ_{1}$, where the subscript refers to the number of
zero-probability cells, corresponding to the tables:

\begin{center}
\begin{tabular}{l|cc|cc}
\multicolumn{1}{c}{} & \multicolumn{2}{c}{$SZ_0$} &
 \multicolumn{2}{c}{$SZ_1$}\\
\emph{Type of cancer} &Female & Male & Female & Male \\ \hline
Lung & + & + & + & +\\
Melanoma & + & + & + & +  \\
Ovarian & + & 0 & + & 0  \\
Prostate & 0 & + & 0 & + \\
Stomach & + & + & 0 & + \\
\end{tabular}
\end{center}


Similarly, for the given data, it is not difficult to realize that there exists  only
one logically consistent instance of the quasi-independence model, i.e. that  having all positive cell-probabilities (except
for the two cells corresponding to structural zeros), which we label
$QI_0$ and is schematically equivalent to $SZ_{0}$ above.

\textbf{4.1. Conventional approach}
We test the model of quasi-independence
against the structural-zero model using a Bayesian approach.
In a ``conventional setting'', wherein no particular
provision for zero-probability cells is envisaged, we would simply
consider one instance for each of the above two models, namely $SZ_0$ and $QI_0$.

Given the cell counts $n = (n_{ij})$, a typical analysis would involve the computation of  the Bayes
factor, see Kass and Raftery (1995),  of $QI_0$ \emph{versus} $SZ_0$, i.e.
\begin{equation}
\BF(QI_0:SZ_0)= \frac {\int
f_{QI_0}(n|\theta_{QI_0})\pi_{QI_0}(\theta_{QI_0})d\theta_{QI_0} }{\int
f_{SZ_0}(n|\theta_{SZ_0})\pi_{SZ_0}(\theta_{SZ_0})d\theta_{SZ_0}}=
\frac{m_{QI_0}(n)}{m_{SZ_0}(n)},
\end{equation}
where
\begin{itemize}
\item
 $f_{SZ_0}$ is the  multinomial sampling distribution
under $SZ_0$, with cell-probabilities $\theta_{SZ_0}=( \teij )$,
$(i,j) \in A$, and
similarly for $f_{QI_0}$  under the quasi-independence model, whose
cell-probabilities  are denoted by $\theta_{QI_0}$; \item $\pi_{SZ_0}$
and $\pi_{QI_0}$ are the prior densities for $\theta_{SZ_0}$,
respectively $\theta_{QI_0}$; \item $m_{SZ_0}$ denote the marginal
distribution of $n$ under $SZ_0$, and similarly for $m_{QI_0}$.
\end{itemize}
To obtain the posterior probability  of model $QI_0$ one should
provide, in addition, its prior probability $p_{QI_0}=\text{Pr}(QI_0)$,
leading to
\begin{equation}
\text{Pr}(QI_0|n)=  \frac{p_{QI_0}\BF(QI_0:SZ_0)}{p_{QI_0}\BF(QI_0:SZ_0)+p_{SZ_0}},
\end{equation}
where $p_{SZ_0}=\text{Pr}(SZ_0)=1-p_{QI_0}$.

A Bayesian analysis of this problem would take the prior $\pi_{SZ_0}$ to be Dirichlet, i.e.
\begin{equation}
\theta_{SZ_0} \sim \Di(\alpha),
\end{equation}
with $\alpha=(\alpha_{ij})$ and $ \alpha_{ij}>0$, see e.g. Bernardo
and Smith (1994, p. 134-5 and 441)  and O'Hagan and Forster (2004,
chapter 12). As a consequence,
$m_{SZ_0}(n)$ is a Multinomial-Dirichlet with distribution

\begin{multline}
\label{eq:multinomial dirichlet}
m_{SZ_0}(n) =
\frac{N!}{\prod_{(i,j)\in A}  n_{ij}!}
\times
\frac{H_A(\alpha)}{H_A(\alpha^*)}, \\
  n=(n_{ij}),  \quad n_{ij}=0,1,\ldots,N, \quad
\sum_{i,j} n_{ij}=N
\end{multline}
where
\begin{equation*}
H_T(y)=\frac{\Gamma(\sum_{ t \in T} y_t)}
{ \prod_{t \in T}\Gamma(y_t)},
\end{equation*}
and $\alpha^*=\alpha+n$.

 Consider now the  quasi-independence model $QI_0$, and in particular the choice of the prior
 $\pi_{QI_0}$. This presents some conceptual and practical challenges, that
 we now try to elucidate. Although, in principle, priors under
 distinct models need not be related, as they express prior beliefs
 conditionally on different states of information, it is nevertheless
 desirable that they should be related at least when models are nested
 within an encompassing model.  Pragmatically, this would simplify the
 elicitation task, since one would only assign a prior on the
 parameter under the latter model, and then derive the corresponding
 priors under each of the remaining models from this single
 prior. This procedure should also achieve some sort of internal
 \lq \lq compatibility\rq \rq{} among prior specifications. A general discussion of
 strategies for building compatible priors under several related
 models is contained in Dawid and Lauritzen (2001).  Further
 discussion, elaboration and references may be found in Consonni et al. (2005), and in Consonni and
 Veronese (2006).

Before turning to model $QI_0$,  it is expedient to
rewrite the joint distribution of the counts $n_{ij},\, (i,j) \in A,$ for the $SZ_0$-model as
\ben
\label{eq:cut}
f_{SZ_0}(n|\theta)=f_{SZ_0,1}(n_{(1)}|\theta)
\times f_{SZ_0,2}(n_{(2)}|n_{(1)},\theta),
\een
where
\be
n_{(1)}&=&(n_{31},n_{42}, N-n_{31}-n_{42})\\
n_{(2)}&=& (n_{ij}: (i,j) \in  A \setminus \{ (3,1),(4,2) \})
\ee
Since,  for $(i,j) \in A$,  the joint distribution of $n=(n_{ij})$, under $SZ_0$,
is multinomial with size $N$ and vector of probabilities
$\theta=(\theta_{ij})$,  written $\mbox{Mu}(N;\theta)$, it is easy to check that
 $f_{SZ_0,1}(n_{(1)}|\theta)$ is a $\mbox{Mu}(N;\lambda)$ with
 \be
\label{eq:reparametrization}
\lambda_1=\theta_{31},\,\lambda_2=\theta_{42},\,\lambda_3=1-\lambda_1-\lambda_2,
\ee
while $f_{SZ_0,2}(n_{(2)}|n_{(1)},\theta)$ is given by $\mbox{Mu}(N-n_{31}-n_{42}; \gamma)$,
where
\be
\label{eq:reparametrization1}
\gamma_{ij}=\frac{\theta_{ij}}{1-\theta_{31}-\theta_{42}}
=\frac{\theta_{ij}}{\sum_{(i,j) \in  A \setminus \{ ((3,1),(4,2)) \}}},
\quad  \,(i,j) \in  A \setminus \{ ((3,1),(4,2)) \}.
\ee
The parameters $\lambda$ and $\gamma$ are variation independent, i.e. their joint range
is the  product of the two individual ranges.

Under $QI_0$ we must have
\be
\gamma_{ij}=\gamma_{i+}\gamma_{+j}, \quad (i,j) \in A \setminus \{ (3,1),(4,2) \},
\ee
where
\be
\gamma_{i+}= \gamma_{i1}+\gamma_{i2}, \quad i=1,2,5\\
\gamma_{+j}=\gamma_{1j}+\gamma_{2j}+\gamma_{5j}, \quad j=1,2.
\ee

Let $\gamma_{R}$ denote the collection of $\gamma_{i+}$, and $\gamma_{C} $ that of $\gamma_{+j}$.
Then the distribution of the counts $n$ under $QI_0$ can be written as
\ben
\label{eq:fQI0}
f_{QI_0}(n|\lambda,\gamma_{R},\gamma_{C})=
f_{QI_0,1}(n_{(1)}|\lambda)f_{QI_0,2}(n_{(2)}|n_{(1)};\gamma_{R},\gamma_{C}),
\een
where $f_{QI_0,1}(n_{(1)}|\lambda)$ is $\mbox{Mu}(N;\lambda_1,\lambda_2,\lambda_3)$ and so coincides with
the expression of $f_{SZ_0,1}(n_{(1)}|\theta)$
 in
\eqref{eq:cut}, while $f_{QI_0,2}(n_{(2)}|n_{(1)},\gamma_{R},\gamma_{C})$ is
given by
\ben
\label{eq:fQI02}
f_{QI_0,2}(n_{(2)}|n_{(1)};\gamma_{R},\gamma_{C})=
\frac{(N-n_{31}-n_{42})!}{\prod_{(i,j)\in A \setminus \{ (3,1),(4,2) n_{ij}! \} }}
\times\gamma_{1+}^{n_{1+}}\gamma_{2+}^{n_{2+}} \gamma_{5+}^{n_{5+}}
\times \gamma_{+1}^{\tilde{n}_{+1}} \gamma_{+2}^{\tilde{n}_{+2}},
\een
where $\tilde{n}_{+j}=n_{1j}+n_{2j}+n_{5j}$.

One can thus see that under $QI_0$ the joint distribution factors into three terms, one involving
$\lambda$, one involving $\gamma_{R}$ and one involving $\gamma_{C}$.

Consider now the prior distribution. Given that $\theta_{SZ_0} \sim \Di(\alpha_{SZ_0})$ we first remark  that $\lambda$ and $\gamma$,
are independent, because of ii) of Lemma \ref{lemma:Dirichlet}, see Appendix;
as a consequence we also get that $\lambda$ is independent of the pair $(\gamma_{R},\gamma_{C})$.
Furthermore $\gamma \sim Di(\alpha_{ij},\, (i,j) \in A \setminus \{ (3,1),(4,2) \})$, so that
  $\gamma_{R} \sim \Di(\alpha_{R}) $ and $\gamma_{C} \sim \Di(\alpha_{C}) $, where
$\alpha_R$ and $\alpha_C$ are defined in accordance with $\gamma_R$ and $\gamma_C$, respectively.
Assuming
independence of $\gamma_{R}$ and $\gamma_{C}$ makes the computation of the marginal distribution
$m_{QI_0}(n)$ straightforward since we can integrate separately  the three terms in
\eqref{eq:fQI0}, see also \eqref{eq:fQI02}, each integral being, up to the multinomial coefficient,
of type Multinomial-Dirichlet.

Specifically we get
\ben
\label{eq:mQI0}
m_{QI_0}(n)&=&\frac{N!}{\prod_{(i,j)\in A}  n_{ij}! } \nonumber\\
 &\times&
 \frac{H(\alpha_{31},\alpha_{42},\alpha_+ - \alpha_{31}-\alpha_{42})}
 {H(\alpha^*_{31},\alpha^*_{42},\alpha^*_+ - \alpha^*_{31}-\alpha^*_{42})} \nonumber \\
 &\times&
 \frac{H(\alpha_{1+},\alpha_{2+},\alpha_{5+})}
 {H(\alpha^*_{1+},\alpha^*_{2+},\alpha^*_{5+})}
 \times
 \frac{H(\tilde{\alpha}_{+1},\tilde{\alpha}_{+2})}
 {H(\tilde {\alpha}^*_{+1},\tilde{\alpha}^*_{+2})},
\een
where
$\alpha_+=\sum_{(i,j)\in A}  \alpha_{ij}$, $\tilde{\alpha}_{+j}=\alpha_{1j}+\alpha_{2j}+\alpha_{5j}$,
$\tilde{\alpha}^*_{+j}=\alpha_{1j}+n_{1j}+\alpha_{2j}+n_{2j}+\alpha_{5j}+n_{5j}$.

\textbf{4.2. Allowing for zero-probability cells}
Phylosophically, we earnestly take the view that each
instance of a model must be assigned \emph{a-priori} a positive
probability: in this sense we completely adhere to the principle that
Lindley (1985, p.104) names ``Cromwell's rule''. This leads us naturally to the
idea of regarding a model $\mathcal M$ as a finite
\emph{mixture} of its instances. This aspect represents a characterizing
feature of our approach to the analysis of contingency tables.

We can thus write the
mixture representation of $\mathcal M$ as
 \begin{equation}
\label{formula:mixture model U} f_{\mathcal M}(n|\theta_{\mathcal
  M})= \sum_h
q_{{\mathcal M}_{h}}f_{{\mathcal M}_{h}}(n|\theta_{{\mathcal M}_{h}}),
\end{equation}
where $\theta_{\mathcal M}$ is the
collection of all instance-specific parameters $\theta_{{\mathcal M}_{h}}$ and
$q_{{\mathcal M}_{h}}$ is the prior probability attached to instance
$\mathcal M_h$.

Specializing \eqref{formula:mixture model U} to the $SZ$ and $QI$ model, and then computing the marginal distribution of the data under each model, leads to the Bayes factor
\begin{equation}
\label{eq:BF QI:SZ}
\BF(QI:SZ)=\frac{ q_{{QI}_{0}} m_{{QI}_{0}}(n)}{q_{{SZ}_{0}}
m_{SZ_{0}}(n)+q_{{SZ}_{1}} m_{SZ_{1}}(n)}.
\end{equation}

Let us now consider in detail the computations that are needed for the
evaluation of $\BF(QI:SZ)$. Let $\xi \in (0,1)$ be the chance that a
cell has zero probability, and assume that the allocation of
zero probability to each cell takes place independently. Then, we can
derive $q_{{SZ}_{0}}$ and $q_{{SZ}_{1}}$, and obtain
\begin{gather}
q_{{SZ}_{0}} = \frac{(1-\xi)^8}{1-\xi^8},\\
q_{{SZ}_{1}} = \frac{\xi(1-\xi)^7}{1-\xi^8}.
\end{gather}
Consider now the assignment of $q_{{QI}_{0}}$. We recall that we have 87
instances with total probability $C(\xi)$, then
\begin{equation}
\label{eq:qQI0}
q_{{QI}_{0}}=\frac{(1-\xi)^8}{C(\xi)}.
\end{equation}

Table \ref{tab:prior model probabilities} reports the value of $q_{{SZ}_{0}}, q_{{SZ}_{1}}, q_{{QI}_{0}}$
 for selected choices of $\xi$
(for values of $\xi$ above 0.5, the values are zero to two decimal places).
\begin{table}
\begin{center}
 \caption{Prior probabilities $q_{{SZ}_{0}}$, $q_{{SZ}_{1}}$, $q_{{QI}_{0}}$ for selected values of $\xi$}
 \label{tab:prior model probabilities}
\begin{tabular}
{|l||l|l|l|}\hline
$\xi$ & $q_{{SZ}_{0}}$ & $q_{{SZ}_{1}}$ & $q_{{QI}_{0}}$\\
\hline
0.1 &  0.43     &0.05 &0.78\\
0.2 & 0.17& 0.04& 0.51\\
0.3 & 0.06& 0.03&0.23\\
0.4 &0.02 & 0.01& 0.07\\
0.5 &0.00 &0.00 &0.01\\
\hline
\end{tabular}
\end{center}
\end{table}

We now consider the marginal distribution of the data under the $SZ_1$-instance.
The conditioning method  of Lemma
\ref{lemma:Dirichlet}, item ii), leads immediately to conclude that
$\theta_{SZ_1} \sim \Di(\alpha_{SZ_1})$, where $\alpha_{SZ_1}= (\alpha_{ij}, \, (i,j) \in A \setminus \{ (5,1)\})$,
whence $m_{SZ_1}$ has an expression analogous to that of $m_{SZ_0}$, the only difference being that now the set
over which the  indexes vary is $A \setminus \{ (5,1)  \}$.

For given $\xi$ and $\alpha$, the Bayes factor $\text{BF}(QI:SZ)$ can now be computed using \eqref{eq:BF QI:SZ}.
Notice that the multiplicative term $\frac{N!}{\prod_{(ij) \in A }n_{ij}!}$ appears
both in the numerator and denominator of \eqref{eq:BF QI:SZ}, and so cancels out (strictly speaking  the
product for the instance $SZ_1$ is over a set that does not contain $(5,1)$: however since
$n_{51}=0$ the result is the same whether this value appears or not).

Consider first the assignment of $\xi$, which represents the chance that a cell has probability zero.
Save for the case of a structural zero, it seems reasonable that we should
assign a low value to $\xi$, since the corresponding event  should be regarded \emph{a priori}
as a rather unusual circumstance.
 In view of Table \ref{tab:prior model probabilities}, setting $\xi=0.1$ seems a  sensible choice. Indeed, while
 the prior probability of model $QI$ is  higher than that of $SZ$, nevertheless the discrepancy between the
 two values (0.78 against 0.48) is less pronounced for this choice of $\xi$ than for other choices, so that
 the comparison between the two models is fairer.

We now take into consideration the choice of $\alpha$.
Unless
there exists substantive prior information  allowing to discriminate
\emph{a priori} between cells,
we shall  choose the same value $\bar{\alpha}$ for each $\alpha_{ij}$; also
low values of $\bar{\alpha}$ are typically  recommended, whenever prior information is weak. Natural choices are
represented by $\bar{\alpha}=0.5$, corresponding to Jeffreys prior, or $\bar{\alpha}=1$, corresponding to a uniform
prior on the simplex.

We now provide a method for the choice of $\bar{\alpha}$, using the technique of the  \emph{imaginary training sample}.
This method has been implemented for instance by
Spiegelhalter and Smith (1980) to deal with model choice using
 improper priors. We believe however that the idea can be usefully applied also in the context of
 proper priors,  see Consonni et al. (2005) for a similar elaboration.

Consider for simplicity only the models $SZ_0$ and $QI_0$.
Suppose we can identify a \emph{minimal imaginary training sample}
that provides \emph{maximal} support (irrespective of the prior) to model $QI_0$.
Then it is reasonable to require that the Bayes factor for these fictitious data should be approximately
1, i.e. the models are \lq \lq equally likely\rq \rq{} in terms of the empirical evidence.
To see why this should be the case, notice that, on the one hand the data actually support $QI_0$ very strongly;
on the other hand,  the sample size is so
small that the evidence in favor of either model should be roughly the same. The condition that
the  Bayes factor should be equal to 1
can be employed to select reasonable values for the hyper-parameters of the prior distribution.

Consider the situation in which we have 1 observation in each cell, for a
total of 8 observations. It is straightforward to  verify that this table is perfectly consistent with
the $QI_0$-model: in particular  the actual and fitted counts (the latter based on ML  estimates)
coincide. If we fix $\xi=0.1$ as suggested above, the value $\bar{\alpha}=1$ provides a Bayes factor equal to 1.03,
which is quite satisfactory;
on the other hand  $\bar{\alpha}=0.5$ would give a BF equal to 0.67. We also experimented with other values of
$\bar{\alpha}$ and did not get values of BF close to 1.

Having set $\xi=0.1$ and  $\bar{\alpha}=1$,
we  now proceed to the analysis  of the cancer data.  The Bayes factor of $QI$ against
$SZ$ is equal to 0.17, which is clearly not supporting the hypothesis of quasi-independence. To better
assess this value, it is useful to derive the Bayes factor \emph{against} $QI$, which is
merely the reciprocal of the above, and to further transform it using the logarithm in base 10. In this way
we can make use of
the scale developed by Jeffreys, see Kass and Raftery (1995) and Robert (2001, p. 228), for the
interpretation of the evidence provided by a Bayes factor.
Specifically, the evidence \emph{against} $QI$ is
\begin{itemize}
\item
\emph{poor} if $0<\log_{10} \BF(SZ:QI)<0.5$,
\item
\emph{substantial} if $0.5<\log_{10} \BF(SZ:QI)<1$,
\item
\emph{strong} if $1<\log_{10} \BF(SZ:QI)<2$,
\item
\emph{decisive} if $\log_{10} \BF(SZ:QI)>2$,
\end{itemize}
where $\BF(SZ:QI)=1/\BF(QI:SZ)$.
As a consequence we get $\log_{10}(1/0.17)=0.77$ which thus represents \emph{substantial} evidence against $QI$,
essentially in accord with the frequentist answer which states a p-value between $2\%$ and $3\%$.
It is instructive
to verify what would have been the result of a conventional Bayesian analysis,
based exclusively on the positive-cell models $SZ_0$ and $QI_0$,
as opposed to the model based on mixtures developed in this paper.
 Recall that, in the standard case, the BF would  simply be  the ratio
$m_{QI_0}(n)/m_{SZ_0}(n)$. In this case the BF takes the value 0.55,
which is appreciably higher than the value 0.17 obtained with
our analysis. More interestingly, when translated to the Jeffreys scale, we obtain
$\log_{10}(1/0.55)=0.26$ which only represents \emph{poor} evidence against $QI$, which
is an order of magnitude lower, on the Jeffreys scale, than the one we obtained with our analysis.

\setcounter{chapter}{4} 
\setcounter{equation}{0} 
\noindent {\bf 4. Discussion}

In this paper we have presented a new methodology for the Bayesian analysis of contingency
tables that allows explicitly for the possibility of zero-probability cells.

The essential features of our approach are: the notion of extended log-linear model,
the support of computational
algebraic geometry to enumerate and list all model-instances having varying support, the use
of a mixture model to represent the sampling distribution of the cell-counts,
a technique to assign prior probabilities to the various model-instances, a
method to derive prior distributions on the parameter space of each model-instance
starting from a Dirichlet prior under the structural zero model, as well as an elicitation procedure for the
corresponding hyper-parameters.

We have illustrated our methodology by means of an application to a real data set
involving a cross-classification of types of cancer and gender. The corresponding contingency table
presents two structural zeros, and a cell with a zero count.
 The results we obtain, when testing the hypothesis of quasi independence, show
 that our methods can  lead to conclusions that are substantively different from
those based on a standard modeling analysis, which does not explicitly allow for the possibility
 of  zero-probability
cells.

In order to apply the  algebraic Bayesian approach presented in this
paper to large and sparse contingency tables, we believe that a
purely ``automated'' approach can be expected to run
into serious computational issues, although technology is rapidly
evolving in this area as for instance evidenced, within the field of Maximum Likelihood Estimation,
in the recent paper by Erikkson et al. (2006), see also Patcher and Sturmfels (2005) for
a variety of high-dimensional applications. A careful choice of prior distribution  is often
the only  sensible way  to make the analysis viable, see for
instance Diaconis and Rolles (2006) in the context of Markov chains
with forced zeros. We therefore believe that a blend of
computational algebraic methods and prior information on the set of
possibly-zero probability cells is likely to be the best option for
the analysis of moderate to large multi-way tables.

\noindent {\bf Appendix}

We summarize below some useful facts about
 the Dirichlet distribution, see e.g. Bernardo and Smith (1994,
 pp. 134-5) (notice however that our notation is slightly different
 from theirs).
\begin{lem}
\label{lemma:Dirichlet}
Let $\theta=(\theta_1,\ldots,\theta_s)$, with $0 < \theta_k < 1$, $k=1,\ldots,s$, and $\sum_{k=1}^s
\theta_k=1$. Assume that $\theta \sim \Di(\alpha)$, with $\alpha=(\alpha_1, \ldots, \alpha_s)$ and $\alpha_k >0$.
\begin{itemize}
\item [i)]
  \begin{equation*}
\left(\theta_1,\ldots, \theta_r, (1-\sum_{l=r+1}^s\theta_l)\right)
\sim \Di \left(\alpha_1,\ldots,\alpha_r, \sum_{l=r+1}^s \alpha_l \right),
\quad r<s.
  \end{equation*}
\item[ii)]
Let $\theta_m^{\prime}=\frac{\theta_m}{\sum_{q=1}^r \theta_q}, \, m=1,\ldots, r, \,
r < s$, then
$$(\theta_1^{\prime}, \ldots, \theta_r^{\prime})
 \sim \Di(\alpha_1, \ldots, \alpha_r ),$$
and $(\theta_1^{\prime}, \ldots, \theta_r^{\prime})$
 is independent of $(\theta_{r+1}, \ldots, \theta_s)$.
 \item[iii)]
 Let $\theta^*_1=\theta_1+\ldots+\theta_{i_1}$, $\ldots$,
 $\theta^*_t=\theta_{i_{t-1}}+\ldots+\theta_{s},\quad 1 \leq t < s$, then
 $$
(\theta^*_1,\ldots, \theta^*_t) \sim \Di(\alpha^*_1,\ldots,\alpha^*_t),
 $$
 $\alpha^*_1=\alpha_1+\ldots+\alpha_{i_1}$, $\ldots$, $\alpha^*_t=\alpha_{i_{t-1}}+\ldots+\alpha_{s}$.
\end{itemize}
\end{lem}

\noindent{\large\bf Acknowledgment}

Work partially supported by MIUR, Rome, under the projects PRIN 2003138887 and PRIN 2005132307, by the
University of Pavia, the University of Genova and Politecnico of Torino. We thank Simplice Dossou-Gb\'et\'e
 and
Laboratoire de Math\'ematiques Appliqu\'es UMR CNRS 5142 at
Universit\'e de Pau et des Pays de l'Adour
  for providing hospitality and support
while part of this article was written. The second author especially
thanks H.P. Wynn for many discussions and suggestions. A special thank to Persi
Diaconis who provided us with thoughtful feedback. Finally, the
careful reading and comments by two referees are  gratefully
acknowledged.

\vspace{.5cm}
\noindent{\large\bf References}

\begin{description}

\item {4ti2 team}. 4ti2 -- A software package for algebraic, geometric and
              combinatorial problems on linear spaces. \verb+http://www.4ti2.de+.

\item Bernardo, J. M. and Smith, A.F.M. (1994). \emph{Bayesian Theory}. Wiley,
Chichester.

\item Bishop, Y. M. M., Fienberg, S. E. and Holland, P. W. (1975). \emph{Discrete Multivariate Analysis}.
MIT Press. Cambridge, MA.

\item Brown, L. D. (1986). \emph{Foundations of Exponential Families}. IMS Lecture Notes-Monograph Series
\textbf{6}. Hayward, CA.

\item Christensen, R. (1997). \emph{Log-linear Models and Logistic
Regression}. Springer, New York.

\item {CoCoA}Team, {{\hbox{\rm C\kern-.13em o\kern-.07em C\kern-.13em
      o\kern-.15em A}}}: a system for doing {C}omputations in
{C}ommutative {A}lgebra. \verb+http://cocoa.dima.unige.it+.

\item Consonni, G. and Veronese, P. (2006). Prior specifications for the
comparison of linear models. Submitted.

\item Consonni, G., Guti\'{e}rrez-Pe\~{n}a, E. and Veronese, P. (2005).
Compatible priors for Bayesian model comparison with an
  application to the Hardy-Weinberg equilibrium model. Under revision for
 \emph{Test}.

\item Dawid, A. P. and Lauritzen, S. L. (2001). Compatible prior
distributions. In \emph{Bayesian Methods with Applications to
  Science, Policy and Official Statistics} (E. George, ed.).
Monographs of Official Statistics, pp.\ 109--118. Office for official
publications of the European Communities: Luxembourg. \\
\verb+http://www.stat.cmu.edu/ISBA/index.html+.

\item Diaconis, P. and Rolles, S. W. W. (2006). Bayesian analysis for
reversible Markov chains. \emph {The Annals of Statistics},
\textbf{34}, 1270-1292.

\item Eriksson, N., Fienberg, S. E., Rinaldo, A., Sullivant, S. (2006). Polyhedral conditions for the
nonexistence of the MLE for hierarchical log-linear models.
\emph{Journal of Symbolic Computation}, \textbf{41}, 222-233.

\item Garcia, Stillman and Sturmfels (2005). Algebraic geometry of
Bayesian networks. \emph{Journal of Symbolic Computation}, \textbf{39}, 331-355.

\item Geiger, D., Heckerman, D., King, H. and Meek,
Ch. (2001). {Stratified exponential families: graphical models and
  model selection}, \emph{The Annals of Statistics}, \textbf{29}, 505-529.

\item Geiger, D, Meek, C. and Sturmefels, B. (2006). On the toric algebra
of graphical models. \emph{The Annals of Statistics}, \textbf{34},
1463-1492.

\item Haberman, S. J. (1974). \emph{The Analysis of Frequency Data}. University of Chicago Press.
Chicago.

\item Kass, R. E and Raftery, A. E (1995). Bayes factors.
\emph{Journal of the American Statistical Association}, \textbf{ 90},
773--795.

\item Lauritzen, S. L.(1975). General exponential models for discrete observations.
\emph{Scandinavian Journal of  Statistics}, \textbf{2}, 23--33.

\item Lauritzen, S. L. (1996). \emph{Graphical Models}. The Clarendon
Press Oxford University Press, New York.

\item Leonard, T. and Hsu J. S. J. (1999). \emph{Bayesian Methods}.
Cambridge University Press, Cambridge.

\item Lindley, D.V. (1985). \emph{Making Decisions} 2nd Ed. John Wiley
\& Sons, London.

\item O'Hagan A. and Forster, J. (2004). \emph{Kendall's Advanced
Theory of Statistics}, \textbf{2B}. \emph{Bayesian Inference}. 2nd edition.  Arnold, London.

\item Patcher, L. and Sturmfels, B. (2005). \emph{Algebraic Statistics
  for Computational Biology}. Cambridge University Press, Cambridge, UK.

\item Pistone, G., Riccomagno, E, Wynn, H. P. (2001) \emph{Algebraic
  Statistics: Computational Commutative Algebra in
  Statistics}. Chapman\&Hall, London.

\item Rapallo, F. (2006). Toric statistical models: parametric and
binomial representations. \emph{Annals of the Institute
  of Statistical Mathematics}. Springer DOI 10.1007/s10463-006-0079-z.

\item  Robert, C. P. (2001)
\emph{The  Bayesian Choice}. 2nd edition. Springer, New York.

\item Simonoff, J. S. (2003). \emph{Analyzing Categorical Data}. Springer,
New York

\item Spiegelhalter, D. J.  and Smith, A. F. M. (1980). Bayes factors
and choice criteria for linear models.  \emph{Journal of the Royal
Statistical Society B}, \textbf{ 42}, 215--220.

\item Sturmfels, B. (1996). \emph{Gr\"obner bases and convex
  polytopes}. American Mathematical Society, Providence, RI.

\end{description}

\vskip .65cm
\noindent
University of Pavia, Italy
\vskip 2pt
\noindent
E-mail: guido.consonni@unipv.it
\vskip .65cm
\noindent
Politecnico di Torino, DIMAT, Corso Duca degli Abruzzi 24, 10129
Torino Italy
\vskip 2pt
\noindent
E-mail: giovanni.pistone@polito.it

\end{document}